\documentclass[12pt]{article}
\usepackage{latexsym}
\usepackage{amsmath,amsfonts,amssymb,epsfig}
\usepackage{amsthm,amscd}
\usepackage{textcomp}

\newtheorem{theorem}{Theorem}[section]
\newtheorem{definition}[theorem]{Definition}
\newtheorem{corollary}[theorem]{Corollary}
\newtheorem{proposition}[theorem]{Proposition}
\newtheorem{lemma}[theorem]{Lemma}

\newtheorem{rem}[theorem]{Remark}
\newenvironment{remark}{\begin{rem} \rm}{\end{rem}}
\newtheorem{ide}[theorem]{Idea}

\newtheorem{exa}[theorem]{Example}

\newtheorem{ques}[theorem]{Question}
\newenvironment{question}{\begin{ques} \rm}{\end{ques}}
\newtheorem{pro}[theorem]{Problem}

\newtheorem{spe}[theorem]{Speculation}

\newcommand{\Lie}{{\cal L}}

\font\bbb=msbm10 scaled 1100
\newcommand{\real}{\mbox{\bbb R}}       

\newcommand{\no}{\noindent}

\title{Tight Beltrami fields with symmetry.}
\bigskip
\bigskip\author{R. Komendarczyk\footnote{Department of Mathematics,
University of Pennsylvania;
e-mail: \emph{rako@math.upenn.edu}}}

\begin{document}

\maketitle

\begin{abstract}
\no Let $M$ be a compact orientable Seifered fibered 3-manifold
without a boundary, and $\alpha$ an $S^1$-invariant contact form on
$M$. In a suitable adapted Riemannian metric to $\alpha$, we provide
a bound for the volume $\text{Vol}(M)$ and the curvature, which
implies the universal tightness of the contact structure
$\xi=\ker\alpha$.
\end{abstract}
\begin{center}
\small {\bf keywords}: contact structures, Beltrami fields, curl
eigenfields, adapted metrics, nodal sets, characteristic
hypersurface, dividing sets.
\end{center}

\section{Introduction.}\label{sec:1}
Recall that a \emph{contact structure} $\xi$ on a 3-manifold $M$
is a fully ``nonintegrable'' subbundle of the tangent bundle of $M$.
If $\xi$ is defined by a global 1-form $\alpha$, namely
$\xi=\ker\alpha$, the nonintegrability condition for $\xi$ can be
conveniently expressed as
\begin{align}\label{eq:nonintegrable}
\alpha\wedge d\alpha\neq 0\, .
\end{align}
The contact structure $\xi$ is called \emph{overtwisted}, if there
exists an embedded disk $D\subset M$, such that $T_p D=\xi_p$ along
$\partial D$, $\xi$ is called \emph{tight} if it is not overtwisted.
If all covers of a contact structure are tight, we call it
\emph{universally tight}. A contact structure which is not universally tight 
is either \emph{overtwisted} or \emph{virtually overtwisted} (i.e. its lift to 
a covering space is overtwisted).

As shown in \cite{Eliashberg89}, tight and
overtwisted structures constitute two different types of isotopy
classes among all contact structures in dimension 3. Recall that two
contact structures $\xi_0$ and $\xi_1$ are isotopic if and only if
there exists a homotopy $\xi_t$, $0\leq t\leq 1$, such that each
$\xi_t$ is a contact structure. Clearly, the equivalence up to
isotopy is stronger then the equivalence up to homotopy of plane
fields. Overtwisted structures are rather ``flexible'', and in each
homotopy class of plane fields there exists an overtwisted
representative (c.f. \cite{Eliashberg89}). On the other hand, tight
contact structures are ``rare'', for instance, on manifolds $S^3$,
$\mathbb{R}P^3$, $S^2\times S^1$ there exists the unique, up to
isotopy, tight contact structure. In addition, there exist
$3$-manifolds which admit no tight contact structures (c.f.
\cite{Etnyre01}).

 A more geometric perspective on the contact structures in dimension
 3 has been initiated by Chern and
Hamilton \cite{Hamilton85}. They showed that contact forms
can be equipped with an \emph{adapted} Riemannian metric $g_\alpha$
(see Section \ref{sec:2}). Their main theorem states that an
arbitrary adapted metric can be conformally deformed to a metric of
constant Webster curvature \cite{Hamilton85}. However, the
questions of relations between the Riemannian geometry and the
tight/overtwisted dichotomy, in dimension 3, have not received much
attention in the literature \cite{Blair02}. One may indicate
work in \cite{Belgun03, Geiges97}, on normal CR-structures,
where tightness of Sasakian manifolds is concluded, and also
\cite{Etnyre-Ghrist00, Etnyre-Ghrist00a, Etnyre-Ghrist02b, GK06} where the hydrodynamical
perspectives on contact geometry, related to the Riemannian
geometry, have been studied. The principal motivation behind this paper can be
formulated as follows
\begin{itemize}
\item[] Find conditions on geometric parameters, such 
as volume, curvature and eigenvalues, for a Riemannian metric adapted to a contact structure $\xi$,
that imply tightness of $\xi$. 
\end{itemize}
We achieve this \emph{geometric tightness} for a certain class of $S^1$-invariant contact
structures on Seifered fibered 3-manifolds (theorems in 
Section \ref{sec:6}). This result may be viewed as a translation of Giroux's 
classification of $S^1$-invariant contact structures on $S^1$-bundles, and their quotients,
 into a condition on the volume, curvature and eigenvalues of $M$. In a nutshell, these theorems describe lower bounds for the volume of $M$ in terms of geometric parameters of $M$ and the magnitude $\|\alpha\|$ of a contact form $\alpha$, which defines overtwisted or virtually overtwisted contact structure $\xi$ on $M$.

In Section \ref{sec:7} we present  concluding remarks and 
a slightly different perspective on the problem of geometric tightness, which is motivated by the work of Etnyre and Ghrist in \cite{Etnyre-Ghrist00, Etnyre-Ghrist00a, Etnyre-Ghrist02b}. 
We also indicate geometric conditions that imply tightness of certain contact 
structures on various products and circle bundles, and 
 obtain the well known result \cite{Kanda97} about universal
tightness of $\xi_n=\ker\{\cos(n z) dx+\sin(n z) dy\}$, $n\in
\mathbb{Z}$, on $T^3$. 

\section{Contact structures and adapted metrics.}\label{sec:2}
In this section we show that, given an $S^1$-invariant contact form
$\alpha$, we may adapt a suitable Riemannian metric to $\alpha$ with
a Killing vector field tangent to the $S^1$-fibers preserving $\alpha$. First, we
recall basic results about adapted metrics.
\begin{definition}
 We say that the Riemannian metric $g_\alpha$ is \emph{adapted} to $\alpha$, provided
\begin{align}\label{eq:beltrami_form}
\ast\,d\alpha=\mu\,\alpha,\qquad \mu\neq 0, \quad \mu\in
C^\infty(M),
\end{align}
where $\ast$ is the Hodge star operator in $g_\alpha$ (it is shown
\cite{Hamilton85} that one may additionally prescribe
$\|\alpha\|=1$, and $\mu=2$).
\end{definition}
Here and throughout the paper, we work in the smooth category of closed
$3$-dimensional Riemannian manifolds. We denote
by $D\,T$ the covariant derivative of a tensor field $T$, in the
Levi-Civita connection of $g_\alpha$, and $\nabla f$ the gradient
vector field of a scalar function $f$. We often use the notation 
$\langle\,\cdot\,,\,\cdot\,\rangle$ for the inner product $g_\alpha(\,\cdot\,,\,\cdot\,)$.
By $\xi$ we denote an
orientable contact structure defined by a global 1-form $\alpha$,
i.e. the contact form on a Riemannian 3-dimensional manifold
$(M,g_\alpha)$. Every such $\alpha$ admits the unique, transverse to
$\xi$, vector field $X_\alpha$ called the \emph{Reeb-field} of
$\alpha$ satisfying,
\begin{gather}\label{eq:reeb}
 \alpha(X_\alpha)=1,\qquad \iota(X_\alpha)\,d\alpha=0,
\end{gather}
where $\iota(X_\alpha)$ is the contraction of $d\alpha$ by
$X_\alpha$. Notice that $X_\alpha$ defines a projection
$\pi_\alpha:TM\mapsto \xi$ on $\xi=\ker\alpha$ via the formula
\begin{equation}\label{eq:reeb-projection}
\pi_\alpha(X)=X-\alpha(X)\,X_\alpha .
\end{equation}
If $\|\alpha\|\neq 0$, and $\alpha$ satisfies
\eqref{eq:beltrami_form} one quickly verifies that $\xi=\ker\alpha$
defines a contact structure. Indeed, the nonintegrability condition
\eqref{eq:nonintegrable} holds
\begin{align}\label{eq:nonintegrability-derived}
  \alpha\wedge d\alpha=\mu\,\alpha\wedge\ast\,\alpha=\mu\,\|\alpha\|^2 \ast 1\neq 0.
\end{align}

When we allow $\|X_\alpha\|$ to be non constant, one expects more ``flexibility'' in
metrics adapted to $\alpha$. Following 
\cite[Example 3.7 on p. 93]{Bryant91} we may argue that in the class of analytic metrics
Equation \eqref{eq:beltrami_form} can always be locally solved for a
nonvanishing 1-form $\alpha$, and an arbitrary choice of a constant
$\mu$. As a result, the hyperbolic metric can be locally
adapted to a contact structure (recall that all contact
structures are locally equivalent up to a diffeomorphism \cite{Geiges06}). In contrast, if the condition $\|X_\alpha\|=1$ is imposed,
the hyperbolic metric cannot be an adapted metric  (c.f.
\cite{Blair02}).

\begin{remark}\label{rmk:beltrami}
On a manifold equipped with an
adapted Riemannian metric $g_\alpha$, the dual vector field $v$ to
a contact form $\alpha$ satisfies the Euler equations for the
inviscid incompressible fluid flow (see \cite{Etnyre-Ghrist00a}). 
Such solutions of the Euler equations are
known as \emph{Beltrami fields}. Clearly, if $\mu$ is constant, Beltrami 
fields are just the eigenfields of the curl operator $\ast\,d$ (c.f. \cite{GK06}).
\end{remark}

\no The following lemma provides a useful characterization of
adapted metrics.

\begin{lemma}\label{th:adapt_lemma}
 Given a contact form $\alpha$, a local choice of a metric $g_\alpha$ adapted
 to $\alpha$ is equivalent to a choice of a
 local orthonormal frame  $\{e_1,e_2,e_3\}$ satisfying
 \begin{itemize}
 \item[(i)] $e_1=v\, X_\alpha$, where $X_\alpha$ is a Reeb field of
 $\alpha$ and $v$ a positive function,
 \item[(ii)] $\xi=\text{span}\{e_2, e_3\}$.
 \end{itemize}
 (one may also define the associated almost complex structure $J:\xi\mapsto\xi$ on $\xi$ in terms of the frame as follows: $J e_2=-e_3$, $J e_3=e_2$.)
\end{lemma}

\begin{proof}
Given an adapted metric $g_\alpha$, we have the unique dual vector field $X$,
such that $\alpha(\,\cdot\,)=g_\alpha(X,\,\cdot\,)$.
 We define $e_1=X/\|X\|$, and choose an arbitrary frame on $\xi$ satisfying \emph{(ii)}.
 For the dual coframe $\{\eta_i\}$ to $\{e_i\}$, Equation \eqref{eq:beltrami_form}
 implies
 \begin{eqnarray*}
  \iota(X)d\alpha & = & \iota(X)\mu\,\ast\alpha=\iota(e_1)\mu\,\|X\|\ast\eta_1\\
  & = & \iota(e_1)\mu\,\|X\|\,\eta_2\wedge\eta_3=0\ .
 \end{eqnarray*}
By \eqref{eq:reeb} we conclude that $e_1=v\,X_\alpha$, for some function $v\neq 0$.

Conversely, let $\{e_i\}$ be an adapted frame to
$\alpha$  satisfying \emph{(i)} and \emph{(ii)}, and $\{\eta_i\}$ the
coframe. We must show that the metric $g_\alpha=\sum_i \eta^2_i$ is
adapted to $\alpha$. By \emph{(ii)}, $e_1\perp \xi$ thus
$\alpha(\,\cdot\,)=g(X,\,\cdot\,)$ for $X=h\,e_1=h\,v\,X_\alpha$ and
$\eta_1=w\,\alpha$, where $v, h, w$ are positive functions. Relations among $v, h, w$ follow from the identities:
\begin{eqnarray*}
  \alpha(X_\alpha) & = &
g(X,X_\alpha) =  h\,v\,\|X_\alpha\|^2=1,\\
e_1 & = & v\,X_\alpha=\frac{X_\alpha}{\|X_\alpha\|},\\
\eta_1(\,\cdot\,) & = & g(v\,X_\alpha,\,\cdot\,) =  \frac{1}{h}\,g(X,\,\cdot\,)\\
   & = & \frac{1}{h}\alpha(\,\cdot\,).
\end{eqnarray*}
Thus,
\[
v=w=h=\frac{1}{\|X_\alpha\|}\ . 
\]
Let $a,b,c$ be the coefficients of
$d\alpha$ in $\{\eta_i\}$, because
\[
\iota(e_1)d\,\alpha=v\,\iota(X_\alpha)d\,\alpha=0
\]
 we obtain
\begin{eqnarray*}
\iota(e_1)d\alpha & = & \iota(e_1)\left[a\eta_1\wedge\eta_2+b\eta_1\wedge\eta_3+c\eta_2\wedge\eta_3\right]\\
& = & a\eta_2+b\eta_3=0.
\end{eqnarray*}
Thus $a=b=0$, and 
\begin{eqnarray*}
d\,\alpha & = & c\,\eta_2\wedge\eta_3\\
 & = & c\ast\eta_1=c\,v\,\ast\alpha\ .
\end{eqnarray*}
Equation \eqref{eq:beltrami_form} follows by defining $\mu=c\,v$.
Because $\alpha\wedge d\alpha\neq 0$, we conclude that $\mu\neq 0$.
\end{proof}

\begin{lemma}\label{th:adapt_lemma_g}
Let $\{e_i\}$ be the frame defined locally as in Lemma
\ref{th:adapt_lemma}, and $\{\eta_i\}$ the coframe.
We have the following formula for the adapted metric $g_\alpha$:
\begin{align}\label{eq:adapted_metric_formula}
g_\alpha(X,Y)=\sum_i
\eta^2_i(X,Y)=\frac{1}{v^2}\,\alpha(X)\alpha(Y)+\frac{2\,v}{\mu}\,d\,\alpha(X,J\,\pi_\alpha
Y),
\end{align}
for any $X, Y$, where 
\[
\mu=v\,d\alpha(e_2,e_3)=v\,\alpha(\left[e_2, e_3\right]),\quad \text{and}\quad
v=\|X_\alpha\|\ .
\]

\end{lemma}

\begin{proof}
\begin{eqnarray*}
d\alpha(\,\cdot\,,J\,\cdot\,)&=&\frac{\mu}{v}\eta_2\wedge\eta_3(\,\cdot\,,J\,\cdot\,)\\
& = & \frac{\mu}{2\,
v}\left[\eta_2(\,\cdot\,)\otimes\eta_3(J\,\cdot\,)-\eta_3(\,\cdot\,)\otimes\eta_2(J\,\cdot\,)\right]\\
&=&\frac{\mu}{2\,v}\left(
\eta^2_2(\,\cdot\,,\,\cdot\,)+\eta^2_3(\,\cdot\,,\,\cdot\,)\right)
\end{eqnarray*}
(the last equality follows from \emph{(iii)} in Lemma
\ref{th:adapt_lemma}). But $\eta_1=\frac{1}{v}\,\alpha$, and
\begin{eqnarray*}
g(\,\cdot\,,\,\cdot\,) & = & \sum_i
\eta^2_i(\,\cdot\,,\,\cdot\,)\\
 & = & \frac{1}{v^2}\,\alpha^2(\,\cdot\,,\,\cdot\,)+\eta^2_2(\,\cdot\,,\,\cdot\,)+\eta^2_3(\,\cdot\,,\,\cdot\,)
\\
& = &\frac{1}{v^2}\alpha^2(\,\cdot\,,\,\cdot\,)+\frac{2\,
v}{\mu}d\alpha(\,\cdot\,,J\pi_\alpha\,\cdot\,).
\end{eqnarray*}
\end{proof}
\no As a corollary we conclude a global existence of adapted
metrics \cite{Hamilton85}.
\begin{corollary}\label{th:global_adapted}
 Given a contact form $\alpha$, one may always adapt the
 Riemannian metric $g_\alpha$ to $\alpha$, such that Equation
 \eqref{eq:beltrami_form} is satisfied on $(M,g_\alpha)$.
\end{corollary}
\begin{proof}
Indeed, by Formula \eqref{eq:adapted_metric_formula} for
$g_\alpha$, it suffices to
choose a global almost complex structure $J:\xi\mapsto \xi$,
$\xi=\ker \alpha$, and a vector field $e_1=v\,X_\alpha$, where $v$
is a positive function. The only issue is to define $J$ globally,
but this may be achieved via an arbitrary choice of a metric $g_\xi$
on $\xi$, and defining $J$ by the $\frac{\pi}{2}$-rotation in
$(\xi,g_\xi)$.
\end{proof}

\no These results lead to the following,
\begin{proposition}\label{th:seifert-adapted}
 Suppose that $M$ is a Seifert fibered 3-manifold, and $\alpha$ an invariant contact 1-form
 on $M$ (i.e. invariant under the action of a nonsingular vector field $X$ tangent to 
 $S^1$-fibers of $M$). Then, 
 \begin{itemize}
 \item[(iv)] There exists an adapted Riemannian metric $g_\alpha$ to
 $\alpha$, such that $X$ is a Killing vector field in $g_\alpha$.
 \item[(v)] Moreover, there is an adapted metric $g'_\alpha$, conformal to
 $g_\alpha$, such that $X$ is a unit Killing vector field in $g'_\alpha$.
 \end{itemize}
\end{proposition}
\begin{proof}
 Assume that $\alpha$ defines a positive contact structure i.e.
 $\alpha\wedge d\alpha>0$ (if $\alpha$ defines a negative contact structure the proof is analogous).
 Recall the observation from \cite[Proposition 1.3 on p. 336]{Nicolaescu98}
 stating existence of a Riemannian metric $g$ on $M$, such that
 $X$ is the unit Killing vector field for $g$. Because $X$
 preserves $\alpha$ and $\xi$, the flow $\varphi^t_X$ of $X$
 maps $\xi_p$ isometrically to $\xi_{\varphi^t(p)}$. Consequently, $\Lie_X g_\xi=0$, where $g_\xi$
 denotes the restriction of $g$ to $\xi$.
 By positivity of $\alpha$ and $\Lie_X d\alpha=0$, we have a positive function $v$ such that
 \begin{gather*}
  2\,v\,d\alpha(\pi_\alpha\,\cdot\,,J\,\pi_\alpha\,\cdot\,)=g_\xi(\pi_\alpha\,\cdot\,,\pi_\alpha\,\cdot\,),
 \end{gather*}
 where $J$ is a rotation by $\frac{\pi}{2}$ in $g_\xi$, and $\pi_\alpha$ is the projection defined in \eqref{eq:reeb-projection}.
 Notice that $\Lie_X v=0$, therefore we may define $g_\alpha$ by Formula
 \eqref{eq:adapted_metric_formula}. The conclusion \emph{(iv)} now
 follows from the tensor product formula for the Lie derivative. 
 
 In order to prove \emph{(v)}, consider $h^2=g_\alpha(X,X)$. If $h=1$, we are done, if $h\neq 1$ define 
 $g'_\alpha=\frac{1}{h^2}\,g_\alpha$. We verify that $g'_\alpha$ is
 adapted by plugging into \eqref{eq:adapted_metric_formula}:
\begin{eqnarray*}
g'_\alpha(X,Y) & = & \frac{1}{h^2}
g_\alpha(X,Y)\\
& = & \frac{1}{h^2\,v^2}\,\alpha(X)\alpha(Y)+\frac{2\,v\,h}{h^3\,\mu}\,d\,\alpha(X,J\,\pi_\alpha
Y)\ .
\end{eqnarray*}
This calculation confirms that $g'_\alpha$ is adapted with
$\mu'=h^3\,\mu$. Moreover, $X$ is a unit vector field in
$g'_\alpha$ and, because $\Lie_X h=0$, $X$ has to be a Killing vector field
in $g'_\alpha$.
\end{proof}
\begin{question}
 Can we find $g_\alpha$, which admits both a unit Killing
 vector field $X$ tangent to the fibers of $M$ and $\mu$ as a constant function?
\end{question}

\section{Characteristic hypersurface as a nodal set.}\label{sec:3}

Among known techniques of contact topology, which allow us to detect a
contact isotopy type of a contact structure $\xi$, is the technique of convex
surfaces and dividing curves, originally introduced by Giroux in
\cite{Giroux91, Giroux01}. We adapt Giroux's technique, which is crucial in our further
investigation. First, we briefly review its basic notions.
\begin{definition}
Recall that a vector field $X$ on $M$ is called the \emph{contact
vector field} for $\xi$ if and only if its flow preserves the plane
distribution $\xi$. The set of tangencies $\Gamma_X=\{p\in M: X_p\in
\xi_p\}$ of $X$ and $\xi$ is called the \emph{characteristic
surface} of $X$ and is denoted by $\Gamma_X$. An embedded surface
$\Sigma$ in $M$ is called the \emph{convex surface} if and only if
there exists a transverse contact vector field $X$ to $\Sigma$. The
set of curves $\Gamma=\Gamma_X\cap \Sigma$ is called the
\emph{dividing set} on $\Sigma$.
\end{definition}
\no Another way to express the condition for $X$ to be a contact vector field is the following 
equation for the contact form $\alpha$:
\begin{gather*}\label{eq:contact_1}
\Lie_X \alpha=h\,\alpha,\qquad\text{for some}\quad h\in
C^\infty(M)\ .
\end{gather*}
The special case occurs when
$h=0$ and the contact field $X$ also preserves the contact form
$\alpha$, we consider this case in our further investigation. Also, notice that the characteristic surface
$\Gamma_X$ is a zero set of the function $f=\alpha(X)$, i.e. $\Gamma_X=f^{-1}(0)$. This function 
is commonly known as the \emph{contact hamiltonian} (c.f. \cite{Geiges06}).

\no Classification of contact structures on $S^1$-bundles over a surface, has been partially achieved by 
Giroux in \cite{Giroux01}, and completed in full generality by Honda \cite{Honda00, Honda00_1}. As it is presented in the following theorem,
$S^1$-invariant contact structures on $S^1$-bundles are fully characterized by the topology of the dividing set 
$\Gamma$ on the base (i.e. projected $S^1$-invariant characteristic surface $\Gamma_{S^1}$) and the Euler number of the bundle.
\begin{theorem}[\cite{Giroux01}]\label{th:giroux_bundle}
Let $\xi$ be an $S^1$-invariant contact structure on the principal
$S^1$-bundle $\pi:P\to \Sigma$, where $\Sigma$ is an orientable
surface. Let $\Gamma=\pi(\Gamma_{S^1})$ be a projection of the
characteristics surface $\Gamma_{S^1}$ onto $\Sigma$. Denote by
$e(P)$ the Euler number of $P$.
\begin{description}
\item[(a)] If $\xi$ is tight and one of the connected components
 of $\Sigma/\Gamma$ bounds a disc, then $\Gamma$ has to be a single circle and $e(P)$ must satisfy
\begin{gather*}
\left\{%
\begin{array}{ll}
    e(P)>0, & \hbox{if } \Sigma\neq S^2\\
    e(P)\ge 0, & \hbox{if } \Sigma=S^2\ .\\
   \end{array}%
\right.\end{gather*}
\item[(b)] For $\xi$ to be universally tight it is necessary and sufficient that one of the following holds
\begin{itemize}
    \item[\emph{(b.1)}] for $\Sigma\neq S^2$, none of the connected components of $\Sigma/\Gamma$ is a
    disc,
    \item[\emph{(b.2)}] for $\Sigma=S^2$, $e(P)<0$ and $\Gamma=\O$,
    \item[\emph{(b.3)}] for $\Sigma=S^2$, $e(P)\ge 0$ and $\Gamma$ is connected.
   \end{itemize}
\end{description}
\end{theorem}
\begin{theorem}[\cite{Giroux01}]\label{th:giroux_product}
Let $\Sigma$ be a convex surface of nonzero genus in the contact
manifold $(M,\xi)$. Let $X$ be the contact vector field transverse to $\Sigma$, and $\Gamma_\Sigma=\Gamma_X\cap \Sigma$ the dividing set of $\Sigma$. Then $\xi$ is
tight, in a tubular neighborhood of $\Sigma$, if and only if none of the components of $\Sigma/\Gamma$ is a
disc.
\end{theorem}
\no These results indicate that the topology of characteristic
surfaces is an indicator of tightness/overtwistedness both on a
local and global level. In the remainder of this section the goal is
to interpret $\Gamma_X$ in the Riemannian geometric setting of
adapted metrics. The following result provides such a
characterization (compare to \cite[Lemma 2.7]{Komendarczyk_nodal}).

\begin{theorem}\label{th:elliptic_contact}
 Assume that $X$ is a global contact vector field on the Riemannian manifold $(M,g_\alpha)$
 which preserves a contact form $\alpha$ satisfying \eqref{eq:beltrami_form}.
 Let $f=\alpha(X)$, denote by $\{e_1=\frac{X}{\|X\|},e_2,e_3\}$ a local adapted
 orthonormal frame and by $\{\eta_1,\eta_2,\eta_3\}$ the dual coframe.
\no Then, the coefficients of
$\alpha=a_k\,\eta_k=\frac{f}{v}\eta_1+a_2\eta_2+a_3\eta_3$
 locally satisfy the first order system
\begin{align}\label{eq:eqn_contact}
\left\{
\begin{array}{ll}
    D_1 f=0,\\
    D_2 f=-\mu\,v\,a_3,\\
    D_3 f=\mu\,v\,a_2,\\
   \end{array}
\right.\end{align} where $v=\|X\|$. Furthermore, $f$ satisfies globally the subelliptic equation:
\begin{align}\label{eq:gener_lap}
 \Delta_E f -\langle \nabla \text{\rm ln}\,h,\nabla
 f\rangle+\mu(\mathcal{E}-\mu)f=0,
\end{align}
where $\mathcal{E}=(\ast d\,\eta_1)(e_1)$, $h=1/\mu v$, and
$\Delta_E$ is the Laplacian on the subbundle $E=\ker \eta_1$. One
may express Equation \eqref{eq:gener_lap} in terms of the global
Laplacian $\Delta_M$, on $M$, as follows
\begin{align}\label{eq:gener_lap_2}
\Delta_M f +\frac{1}{v^2}\nabla^2 f (X,X) -\langle \nabla\,
\text{\rm ln}\, \Bigl(\frac{1}{\mu v}\Bigr), \nabla
f\rangle+\mu(\mathcal{E}-\mu)f=0.
\end{align}
\end{theorem}

\begin{proof}
 Recall formulas \cite{Jost02} for the Hessian and
the Laplacian in a frame $\{e_i\}$:
\begin{eqnarray}\label{eq:hessian1}
 \nabla^2_{i j} f & = & \nabla^2 f(e_i, e_j) = D_i D_j f + \sum_k D_k
 f\,\omega_{i\,k}^j,\\
 \notag \Delta_M f & = & \text{tr}\, \nabla^2 f=\sum_i \nabla^2 f(e_i,
 e_i),
\end{eqnarray}
and the following definitions (summation is assumed over the repeating indices)
\begin{eqnarray*}
 \label{eq:formula1}D_i e_j & = & \omega^k_{ij} e_k, \quad
 \omega^k_i=D_i\eta_k= -\omega^k_{i j}\,\eta_j,\quad \omega^k_{ij}=-\omega^j_{ik}, \\
 \label{eq:formula2} d\,\alpha & = & \eta_i\wedge D_i\alpha,\\
 \label{eq:formula3}D\alpha & = &da_k\otimes\eta_k+a_kD\eta_k=da_k\otimes\eta_k-a_k\omega^k_j\otimes\eta_j,\\
 \label{eq:formula4}\Delta_M & =& -D_i D_i+\omega^j_{i\,i}D_j,
 \end{eqnarray*}
 where $D_i\equiv D_{e_i}$. The proof of Theorem \ref{th:elliptic_contact} is a calculation in the adapted coframe $\{\eta_i\}$. Using Cartan's formula and Equation \eqref{eq:beltrami_form} we
obtain (for $v=\|X\|$):
\begin{eqnarray*}
 0 & = & \Lie_{X} \alpha = \iota(X)d\,\alpha+d\,f\\
   & = & \mu\,\iota(X)\ast\alpha + D_i f\,\eta_i,
\end{eqnarray*}
and 
\begin{eqnarray*}   
-D_i f\,\eta_i & = & \mu\,v\,\iota(X_1)\ast\alpha\\
-D_1 f\,\eta_1-D_2 f\,\eta_2-D_3 f\,\eta_3 & = &
\mu\,v(-a_2\,\eta_3+a_3\,\eta_2)\ .
 \end{eqnarray*}
\no These expressions lead to the following equations
 \begin{align}\label{eq:eq_m1}
 \left\{%
\begin{array}{ll}
    D_1 f=0, \\
    D_2 f=-\mu\,v\,a_3,\\
    D_3 f=\mu\,v\,a_2,\\
   \end{array}%
\right.
\end{align}
and
\begin{eqnarray*}
d\,\alpha & = & \sum_{i<j}
a_{ij}\eta_i\wedge\eta_j  =  \eta_i\wedge D_i\alpha\\
  & = & \eta_i\wedge(D_i
a_k\eta_k - a_k\omega^k_{ij}\eta_j)\\
 & = & D_i a_k\eta_i\wedge\eta_k -
 a_k\omega^k_{ij}\eta_i\wedge\eta_j\ .
\end{eqnarray*}
\no Collecting terms in front of $\ast\eta_1=\eta_2\wedge\eta_3$, we have 
\begin{eqnarray*}
a_{23} & = & D_{2}a_3-D_{3}a_2+a_k(\omega^k_{32}-\omega^k_{23})\\
& = & D_{2}a_3-D_{3}a_2+a_1(\omega^1_{32}-\omega^1_{23})-a_2\omega^2_{23}+a_3\omega^3_{32}\
.
\end{eqnarray*}
Applying Equations \eqref{eq:beltrami_form} and \eqref{eq:eq_m1}:
\begin{align*}
a_{23} = \frac{\mu}{v}\, f =  -D_{2}(\frac{1}{\mu\,v}\,D_2
f)-D_{3}(\frac{1}{\mu\,v}D_3 f)+\frac{f}{v}(\omega^1_{32}-\omega^1_{23})-\frac{1}{\mu\,v}D_3
f\omega^2_{23}-\frac{1}{\mu\,v}D_2 f\omega^3_{32},
\end{align*}
and distributing terms we obtain (for $h=1/(\mu\,v)$)
\begin{eqnarray*}
\mu^2 h\, f  & = &  h(-D_2D_2\,f-D_3D_3\,f+\omega^3_{22}D_3\,f+\omega^2_{33}D_2\,f)\\
         &  & \qquad+\,\mu\,h\,f(\omega^1_{32}-\omega^1_{23}) - D_2 h\,D_2 f - D_3 h\, D_3 f\, .
\end{eqnarray*}
Dividing the above equation by $h$ yields
\begin{equation}\label{eq:ellip_f}
(\Delta_E + L + \nu)f  =  0,
\end{equation}
where
\begin{eqnarray*}
 \Delta_E & = &-D_2 D_2-D_3 D_3+\omega^3_{22}D_3+\omega^2_{33}D_2,\\
 L & = & -\frac{1}{h}(D_2 h\,D_2+D_3 hD_3)=-\langle \nabla\, \text{\rm ln} h, \nabla\,\cdot\,\rangle\\
 \nu & = & \mu(\omega^1_{32}-\omega^1_{23}-\mu)=\mu(\mathcal{E}-\mu),\\
 \mathcal{E} & = & \iota(e_1)\ast d\eta_1\ .
\end{eqnarray*}
Applying $D_1 f=0$ and $D_1 e_1=\omega^k_{11} e_k$, in Equation
\eqref{eq:hessian1}, we express \eqref{eq:ellip_f} in terms of the
Laplacian $\Delta_M$ on $M$:
\begin{eqnarray*}
\Delta_E f & =& \Delta_M f-\langle \nabla f, D_1
e_1\rangle\\
 & = & \Delta_M f+\nabla^2 f (e_1,e_1)=\Delta_M
f+\frac{1}{v^2}\nabla^2 f (X,X),
\end{eqnarray*}
where in the second equation we noticed that 
\[
 \langle \nabla f, D_1 e_1\rangle+\langle D_1 \nabla f,
 e_1\rangle=D_1 \langle \nabla f, e_1\rangle=0,
\]
and
\[ 
 \nabla^2 f (e_1,e_1)=-\langle \nabla f, D_1 e_1\rangle\ .
\]
\end{proof}
\begin{corollary}
If the contact field $X$ is a unit vector field in the metric and $\mu\equiv\text{const}$, Equation \eqref{eq:gener_lap} becomes
\begin{align}\label{eq:gener_lap_1}
\Delta_E\,f+\mu(\mathcal{E}-\mu)\,f=0.
\end{align}
\end{corollary}

\no In the remaining part of this section we work under the assumptions of Theorem \ref{th:elliptic_contact}.

\begin{theorem}\label{th:char_surface_geom}
 A characteristic surface $\Gamma_X=f^{-1}(0)$ is the zero set
(also known as the nodal set) of the solution $f$ to the subelliptic
equation \eqref{eq:gener_lap_2}, and consists of a finite disjoint union of
smooth 2-tori: $\Gamma_X\cong \bigsqcup_i T_i^2$, $T^2_i\cong
S^1\times S^1$.
\end{theorem}
\begin{proof}
 Since $f$ is invariant, under a nonsingular vector field $X$,
regular level sets of $f$ must be 2-dimensional tori.
Clearly, $\Gamma_X$ cannot contain a singular point $p$, since it would imply $\alpha(p)=0$ which
contradicts the contact condition \eqref{eq:nonintegrable}. Now, the
claim follows from Theorem \ref{th:elliptic_contact}.
\end{proof}
\begin{remark}\label{rem:irregular_nodal}
Equation \eqref{eq:gener_lap} can be modified to become an elliptic
equation by adding the term $D_1 D_1 f$. It follows, from the methodology in \cite{Bar97}, that
given any solution $f$ to \eqref{eq:gener_lap}, or
\eqref{eq:gener_lap_1}, the nodal set $N=f^{-1}(0)$ is a union
$N=N_{\text{sing}}\cup N_{\text{reg}}$ of the singular part
$N_{\text{sing}}$ of codimension at least $2$ and the regular part
$N_{\text{reg}}$ which is a codimension $1$ submanifold. We remark
that, in general, the nodal set of a solution to an elliptic partial differential equation
can be very irregular. For example, any closed subset
$A\subset\mathbb{R}^n$ is a nodal set of a
solution to some elliptic equation, \cite{Bar97}.
\end{remark}
\no The following proposition is a
standard result from the elliptic theory \cite{Gilbarg83, Evans98} and Aronszajn's Unique Continuation Principle (see
\cite[p. 235]{Aronszajn57}):
\begin{proposition}\label{th:necc_contact}
 Equation \eqref{eq:gener_lap_2} admits nontrivial solutions, if and only if 
 \begin{align}\label{eq:ineq_existence}
  \mu(\mu-\mathcal{E})\geq 0\qquad \text{on}\ M.
 \end{align}
 Moreover, if the above inequality is strict then $f$ cannot be a locally constant function.
\end{proposition}

\begin{corollary}\label{th:tori_fillable}
If $\mu(\mu-\mathcal{E})>0$ then $M$ is fillable almost everywhere by
$2$-tori.
\end{corollary}

\section{Geometry of the dividing set.}\label{sec:4}
Our strategy for this part is to investigate further
how topology of dividing sets is ``controlled'' by the geometry of the underlying manifold. 
These considerations are essential for the proof of the main theorem, where Equation
\eqref{eq:gener_lap} ``projects'' onto an orientable surface
$\Sigma$ and reduces to the eigenequation
\begin{align}\label{eq:eigen-problem}
 \Delta_\Sigma f=\lambda\,f,\qquad f\in C^\infty(\Sigma),\qquad \lambda\in \real_+.
\end{align}
 We seek conditions on $g_\Sigma$ which imply
that $f^{-1}(0)$ is a homotopically essential collection of curves.
In context of Theorem \ref{th:giroux_bundle} and
\ref{th:giroux_product}, these conditions determine, under appropriate assumptions, tightness of the underlying
contact structure. The next result is inspired by \cite[Lemma 11]{Savo01}.

\begin{proposition}\label{th:nodal_curves_cond}
Let $f$ be a solution to \eqref{eq:eigen-problem} and $\Omega$ a
domain in $\Sigma\setminus f^{-1}(0)$. Denote by $K^\pm$ the
positive (negative) part of the scalar curvature $K$ of
$(\Sigma,g_\Sigma)$. If $\Omega$ is diffeomorphic to a 2-disc $D^2$
with smooth boundary then
\[
4 \pi - 2\int_\Omega K^+\leq \lambda\, \text{\rm Vol}(\Omega).
\]
\end{proposition}
\begin{proof}
The Euler characteristic $\chi(\Sigma)$ of a closed orientable
surface $\Sigma$ can be computed via the celebrated Gauss-Bonnet
formula
\begin{align}\label{eq:gauss-bonnet}
2\pi\chi(\Sigma)=\int_\Sigma K +\int_{\partial \Sigma} \kappa_{\nu}\
,
\end{align}
where $K$ is the scalar curvature, and $\nu$ the unit outward normal
along $\partial \Sigma$. Given a smooth function $f$, on $\Sigma$,
every regular level set $N=f^{-1}(c)$ is a codimension 1 submanifold
in $\Sigma$. In \cite{Savo01} the following formula for the mean
curvature of $N$ has been obtained
\begin{align}\label{eq:level-set}
 H_\nu=\frac{\Delta_N f}{\|\nabla f\|}+\frac{1}{2}
\langle\nabla\text{\rm ln}\|\nabla f\|^2,\nu\rangle,
\end{align}
where $\mathbf{\nu}=\frac{\nabla f}{\|\nabla f\|}$ is pointing
towards $\{f
> c\}$, $\Delta_N$ is the scalar Laplacian on $N$, and $H_\nu$
the mean curvature in the $\nu$ direction (c.f. \cite{Chavel84}).

 By Equation \eqref{eq:eigen-problem} and \eqref{eq:level-set}, we
obtain a formula for the geodesic curvature of $\partial \Omega$:
\begin{eqnarray*}\label{eq:geod_curv_grad}
\kappa_{\nu}=H_\nu & = & \frac{\Delta_\Sigma f}{\|\nabla
f\|}+\frac{1}{2} \langle \nabla\,\text{\rm ln}\|\nabla
f\|^2,\nu\rangle\\
& =& \frac{1}{2} \langle \nabla \,\text{\rm ln}\|\nabla
f\|^2,\nu\rangle,
\end{eqnarray*}
where $\nu=\frac{\nabla f}{\|\nabla f\|}$ points towards $\{f>0\}$,
and because $f\upharpoonright_{\partial \Omega}=0$, the last equality is 
a consequence of \eqref{eq:eigen-problem}. Assume that $f>0$ on
$\Omega$, so that $-\nu=\nu_{\text{out}}$ points outwards (it can be
done without loss of generality since both $f$ and $-f$ satisfy
\eqref{eq:eigen-problem}). Define
\begin{align*}
q=(\|\nabla f\|^2+\frac{\lambda}{2} f^2).
\end{align*}
Clearly, $q\upharpoonright_{\partial \Omega}=\|\nabla f\|^2$ and as a
result we have  $\kappa_\nu=\frac{1}{2} \langle \nabla\,\text{\rm
ln}\,q,\nu\rangle$. Dong's theorem \cite{Dong92} implies the following
estimate for the function $q$:
\begin{align}\label{eq:dongs_estimate}
\Delta\,\text{\rm ln}\,q\leq \lambda-2\,K^{-},\qquad K^{-}=\text{\rm
min}(K, 0).
\end{align}
 By Green's formula (see e.g. \cite[p. 7]{Chavel84}), we obtain
\begin{eqnarray*}
\frac{1}{2}\int_\Omega \text{div}\circ \nabla\,\text{\rm ln}\, q & = & \int_{\partial
\Omega} \frac{1}{2}\langle \nabla \,\text{\rm ln}\,
q,\nu_\text{out}\rangle=-\int_{\partial
\Omega} \frac{1}{2}\langle \nabla \,\text{\rm ln}\, q,\nu\rangle,\\
\frac{1}{2}\int_\Omega \Delta\,\text{\rm ln}\, q & = & \int_{\partial
\Omega} \kappa_{\nu} .
\end{eqnarray*}
(Note that $\Delta=-\text{div}\circ \nabla$, and the orientation
$e_1$ of $\partial \Omega$ is chosen, so that $\{\nu_\text{out},e_1\}$ agrees
with the orientation of $\Omega$.) Applying estimates \eqref{eq:dongs_estimate} and
\eqref{eq:gauss-bonnet}, we derive
\begin{eqnarray*}
\int_{\partial\Omega} \kappa_\nu\ & \leq & \frac{\lambda}{2}\,
\text{Vol}(\Omega)-\int_\Omega K^-\\
2\pi\chi(\Omega)-(\int_\Omega K^+ + \int_\Omega K^-)& \leq &
\frac{\lambda}{2}\, \text{Vol}(\Omega)-\int_\Omega K^-\\
2\pi\chi(\Omega) - \int_\Omega K^+& \leq & \frac{\lambda}{2}\,
\text{Vol}(\Omega).
\end{eqnarray*}
Since $\chi(D^2)=1$, the claim follows from \eqref{eq:gauss-bonnet}.
\end{proof}

\begin{corollary}\label{th:nodal_curves_cond_neg}
Let $f$ satisfy Equation \eqref{eq:eigen-problem}, if $\Sigma$ is a nonpositively curved surface (i.e. $K\leq 0$) of area
$\text{\emph{Vol}}(\Sigma)$, then the following is a necessary condition for
one of the domains in $\Sigma\setminus f^{-1}(0)$ to be a disc with smooth boundary,
\[
\frac{4\pi}{\text{\emph{Vol}}(\Sigma)}\leq \lambda\ .
\]
\end{corollary}

\no In the remaining part of this section, we focus on the case of a convex surface $\Sigma$ embedded in $(M, \xi ,g_\alpha)$,
$\xi=\ker\alpha$. These results are interesting in
their own right, and later provide an essential ingredient in the
proof of the main theorem. 
\begin{proposition}\label{th:ellip_perp}
 Assume the setup of Theorem \ref{th:elliptic_contact}, let $\Sigma$ be a convex surface embedded in $(M,\xi)$, $\xi=\ker\,\alpha$. If a contact field $X$ is
orthogonal to $\Sigma$ and $K\leq 0$, the sufficient condition for
tight tubular neighborhood of $\Sigma$ reads
\begin{align}\label{eq:tight-tub}
\max_{\Sigma} \bigl(\Delta_\Sigma \text{\rm
ln}\|\alpha\|\bigr) < \frac{2\pi}{\emph{\text{Vol}}(\Sigma)}\ .
\end{align}
\end{proposition}
\begin{proof}
Since $X\perp \Sigma$, Equation \eqref{eq:gener_lap} simplifies as
\begin{align}\label{eq:ellip_perp}
\Delta_\Sigma\,f+\frac{\langle \nabla\, \mu\,v, \nabla
f\rangle}{\mu\,v}-\mu^2\,f=0\ .
\end{align}
\no In order to show \eqref{eq:ellip_perp}, we must prove: $\Delta_E=\Delta_\Sigma$ in the frame $\{e_1=\frac{X}{\|X\|},e_2,e_3\}$,
where $\{e_2,e_3\}$ span $T\Sigma$. Equation \eqref{eq:ellip_f}
 yields ($E=T\Sigma$)
\[
\Delta_E=-D_2 D_2-D_3 D_3+\omega^3_{22}D_3+\omega^2_{33}D_2.
\]
Local vector fields $\{e_2,e_3\}$ are tangent to $\Sigma$ and thus the bracket
$[e_2,e_3]$ satisfies: $[e_2,e_3]\in T\Sigma$.
 By the general formula for Christoffel symbols in the frame
\cite{Jost02}:
\begin{align}\label{eq:christoffel_bracket}
\omega^k_{ij}=\frac{1}{2}\{\langle[e_i,e_j],e_k\rangle-\langle[e_j,e_k],e_i\rangle+\langle[e_k,e_i],e_j\rangle\},
\end{align}
we conclude that the formula $\Delta_\Sigma=-D_i
D_i+\omega^j_{i\,i}D_j$ implies $\Delta_E=\Delta_\Sigma$ on
$\Sigma$. In addition,
\begin{eqnarray*}
\langle [e_2,e_3],e_1\rangle & = & \eta_1([e_2,e_3])=0,\\
 d\eta_1(e_2,e_3) & = & 0,\\
\mathcal{E} & = & (\ast d\eta_1)(e_1)=0.
\end{eqnarray*}
Secondly, we express the middle term in \eqref{eq:gener_lap} as follows ($h=1/(\mu v)$):
\begin{align*}
-\langle \nabla\, \text{\rm ln}\, h, \nabla f\rangle=\langle-\mu
v\nabla\Bigl(\frac{1}{\mu v}\Bigr), \nabla
f\rangle=\frac{1}{\mu\,v}\langle \nabla (\mu\,v), \nabla f\rangle,
\end{align*}
which leads to Equation \eqref{eq:ellip_perp}.

In the next step, we calculate the geodesic curvature of $\partial\Omega$, where $\Omega$ is a domain in
$\Sigma\setminus f^{-1}(0)$, and $f=\alpha(X)$ is a solution to Equation
\eqref{eq:ellip_perp}. By \eqref{eq:level-set}
\[
  \kappa_\nu=\frac{\Delta_{\Sigma} f}{\|\nabla f\|}+\frac{1}{2}
  \langle \nabla \,\text{\rm ln}\|\nabla f\|^2,\nu\rangle.
\]
Equations \eqref{eq:ellip_perp} and \eqref{eq:level-set} yield
\begin{align}\label{eq:k_prop1}
  \kappa_\nu=-\frac{\langle \nabla\, (\mu\,v), \nabla
f\rangle}{\mu\,v\|\nabla f\|}+\frac{\mu^2\,f}{\|\nabla
f\|}+\frac{1}{\|\nabla f\|} \langle \nabla\|\nabla f\|,\nu\rangle.
\end{align}
Let $\alpha=a_i\eta_i$, \eqref{eq:eqn_contact} implies that $D_1 f=0$
and
\begin{eqnarray*}
\|\alpha\|^2 & = &\sum_i
a^2_i=\Bigl(\frac{f}{v}\Bigr)^2+\Bigl(\frac{D_2 f}{\mu
v}\Bigr)^2+\Bigl(\frac{D_3 f}{\mu v}\Bigr)^2,\\
(\mu\, v \|\alpha\|)^2 & = & (\mu f)^2+\|\nabla f\|^2,\\
v^2 & = & \frac{f^2}{\|\alpha\|^2}+\frac{\|\nabla f\|^2}{(\mu \|\alpha\|)^2}.
\end{eqnarray*}
Because $f\upharpoonright_{\partial \Omega}=0$, we derive
\begin{eqnarray*}
\mu\,v\upharpoonright_{\partial \Omega} & = &
\frac{\|\nabla\,f\|}{\|\alpha\|},\\
\nabla (\mu v)\,\upharpoonright_{\partial \Omega} & = & -\frac{1}{\|\alpha\|^2}(\nabla\|\alpha\|)\|\nabla f\|+\frac{1}{\|\alpha\|}\nabla\|\nabla f\|,
\end{eqnarray*}
and for $\nu=\frac{\nabla f}{\|\nabla f\|}$:
\begin{eqnarray*}
\frac{\langle \nabla\, (\mu\,v), \nabla f\rangle}{\mu\,v\|\nabla f\|}\Bigl|_{\partial \Omega}
& = & \frac{1}{\mu v}\langle \nabla \mu v, \nu\rangle\\
& = & -\frac{1}{
\|\alpha\|}\langle \nabla \|\alpha\|,\nu\rangle+\frac{\langle
\nabla\|\nabla f\|,\nu\rangle}{\|\nabla f\|}.
\end{eqnarray*}
Substituting in  \eqref{eq:k_prop1} yields
\begin{eqnarray*}\label{eq:geod-curv2}
\kappa_\nu & = & \frac{1}{\|\alpha\|}\langle \nabla
\|\alpha\|,\nu\rangle-\frac{\langle \nabla\|\nabla
f\|,\nu\rangle}{\|\nabla f\|}+\frac{\mu^2\,f}{\|\nabla
f\|}+\frac{\langle \nabla\|\nabla f\|,\nu\rangle}{\|\nabla f\|}\\
 & = &\langle \nabla\,\text{\rm ln}\|\alpha\|,\nu\rangle.
\end{eqnarray*}
\no  As before, we apply \eqref{eq:gauss-bonnet} and the Green's
formula to obtain
\[
 2\pi\chi(\Omega)=\int_\Omega K + \int_\Omega \Delta_\Sigma
 \text{\rm ln}\|\alpha\|\,.
\]
Clearly, if $K\leq 0$ the expression \eqref{eq:tight-tub}, provides a
sufficient condition for $\Sigma$ to have a tight tubular
neighborhood.
\end{proof}
\no Therefore, if we have an
orthogonal contact vector field $X$ to an embedded convex surface $\Sigma$,
Equation \eqref{eq:tight-tub} provides a condition for a tight
tubular neighborhood of $\Sigma$. Contrary to the method of general
convex surfaces, convex surfaces admitting an orthogonal contact
field, as described in Proposition \ref{th:ellip_perp}, are special
(see Remark \ref{rem:convex_special}). In such circumstances, $X$ is
tangent to $\xi$ along the dividing set $\Gamma_X$, and the
orthogonality assumption $X\perp \Sigma$ forces $X_\alpha$ to be
tangent to $\Sigma$, because $X_\alpha\perp \xi$. Based on
Equations \eqref{eq:eqn_contact}, we conclude that the
Reeb field $X_\alpha$ is tangent to the dividing set
$\Gamma_\Sigma$, and $\Gamma_\Sigma$ is a set of periodic orbits of
$X_\alpha$. We have proved, 

\begin{proposition}\label{th:dividing_periodic}
For an embedded surface $\Sigma$, in the contact manifold $(M,\xi)$, satisfying assumptions of Proposition \ref{th:ellip_perp}, the dividing set $\Gamma_\Sigma$ is a set of
periodic orbits of the Reeb field $X_\alpha$.
\end{proposition}

\begin{remark}\label{rem:convex_special}
Example in \cite[p. 327]{Geiges06} demonstrates that the dynamics
of $X_\alpha$ may change drastically depending on a
choice of a contact form $\alpha$ defining $\xi$. Consider the
following family of contact forms on $S^3\subset \real^4$, for $t\ge 0$:
\begin{eqnarray*}
\alpha_t & = & (x_1\,d y_1-y_1\, d x_1)+(1+t)(x_2\,d y_2-y_2\, d
x_2),\\
X_{\alpha_t} & = & (x_1\,\partial_{y_1}-y_1\,\partial_{x_1})+\frac{1}{1+t}(x_2\,\partial_{y_2}-y_2\,\partial_{x_2}),
\end{eqnarray*}
where we consider $S^3$ as a unit sphere in the standard coordinates $(x_1,x_2,x_3,x_4)$ in $\real^4$.

\no If $t=0$, $X_{\alpha_0}$ defines a Hopf fibration on $S^3$, in
particular, all orbits of $X_{\alpha_0}$ are closed. For $t\in
\real\backslash\mathbb{Q}^+$, $X_{\alpha_t}$ defines an irrational
flow on tori of the Hopf fibration and has just two periodic
orbits (at $x_1=y_1=0$, and $x_2=y_2=0$).
It demonstrates that in the irrational case any embedded surface
away from the periodic orbits cannot admit the contact vector
required in Proposition \ref{th:dividing_periodic}. (It also
demonstrates that contact forms are not stable, i.e. in the above
example there exists no family of diffeomorphisms $\psi_t$ such
that $\psi_{t\,\ast} \alpha_t=\alpha_0$, otherwise the flows of
$X_{\alpha_t}$ would have to be conjugate).
\end{remark}

\section{Riemannian submersions and the horizontal Laplacian.}\label{sec:5}

Proposition \ref{th:ellip_perp} describes situations where the
horizontal Laplacian $\Delta_E$, in Equation \eqref{eq:ellip_f}, becomes the
Laplacian on a surface.  We begin by proving a similar
statement in the setting of a Riemannian submersion on a principal $S^1$-bundle. The
reader may consult \cite{Gilkey98} for the general treatment of
related questions for the Hodge Laplacian on forms. 

A submersion $\pi:M\longrightarrow N$ is Riemannian if and only if 
\[
\pi^\ast:T_p
M\supset\ker(\pi^\ast)^\perp_p\longrightarrow T_{\pi(p)} N,
\] 
determines a linear isometry, for all $p\in M$. In other words, for $V, W\in TM$ which are perpendicular
to the kernel of $\pi^\ast$, we have $g_M(V,W)=g_N(\pi^\ast\,
V,\pi^\ast\, V)$. Every Riemannian submersion
determines an orthogonal decomposition $TM=V\oplus H$ of the tangent
bundle into a vertical subbundle $V=\text{Ker}(\pi^\ast)$ and a
horizontal subbundle $H=V^{\perp}$. 
The main feature of $\pi$ is a possibility of
lifting orthogonal frames on $N$ to horizontal vectors on $M$, which
stay mutually orthogonal. Consequently, we may complete a lifted
frame to an orthogonal frame on $M$.

We
summarize useful, for us, properties of Riemannian submersions through a series of lemmas, where vectors on the base $N$ are denoted with capital
letters $E, F$ and lifted vectors on $M$ by small letters $e, f$. We summarize  properties of  the \emph{horizonal lift} operation, $\mathfrak{H}:T_{\pi(p)} N\to T_p M$, in the following (c.f. \cite{Gilkey98}). 

\begin{lemma}[\cite{Gilkey98}]\label{th:bracket_submersion}
Let $\pi: M\to N$ be a Riemannian submersion then
\begin{itemize}
\item[(a)] Lifted $f_p=\mathfrak{H}(F_{\pi(p)})$ is horizontal i.e. $f_p\in
H_p$.
\item[(b)] For any point $p\in M$ and a vector $F_{\pi(p)}\in
T_{\pi(p)} N$, $\pi^\ast \mathfrak{H}(F_{\pi(p)}) =F_{\pi(p)}$.

\item[(c)] Let $f_i=\mathfrak{H}(F_i)$, then $\pi^\ast([f_1,f_2])=[F_1,
F_2]$.

\item[(d)] Let $D^M_i e_j  =  \omega^k_{ij} e_k$, and $D^N_a
E_b  =  \Omega^c_{ab} E_c$. Christoffel symbols satisfy
\begin{align}\label{eq:christ_lift}
\omega^c_{a b}=\Omega^c_{a b}\circ\pi\ .
\end{align}
\end{itemize}
\end{lemma}


\begin{lemma}\label{th:princ_submersion}
Suppose $\pi:P\to \Sigma$ is a projection of an $S^1$-bundle $P$,
equipped with a Riemannian metric $g_P$, which admits a vertical unit
Killing vector field $X$. We have the following,
\begin{itemize}
\item[(e)] $\pi$ defines the Riemannian submersion with an appropriate
choice of the metric on $\Sigma$.
\item[(f)] In a local orthogonal frame of vector fields
$\{e_1,e_2,e_3\}$, where $e_1=X$ and $\{e_2=\mathfrak{H}(E_2), e_3=\mathfrak{H}(E_3)\}$ is the horizontal
lift of a frame $\{E_2, E_3\}$ from $\Sigma$:
\begin{eqnarray}
 \label{eq:princ_sub1} [e_1,e_k] & = & 0,\qquad k=1,2,3,\\
 \label{eq:princ_sub2}\pi\circ\Delta_E & = & \Delta_\Sigma\circ\pi,
\end{eqnarray}
$\Delta_\Sigma$ denotes the Laplacian on $\Sigma$, and $\Delta_E$
is defined in \eqref{eq:ellip_f}, where $E=\text{span}\{e_2,e_3\}$.
\end{itemize}
\end{lemma}
\begin{proof}
Since $X$ is a unit Killing vector field, its flow $\phi^t$ is a
flow of isometries on $P$. Therefore, in a local trivialization:
$(t,\mathbf{x})\in V\cong S^1\times U$, $\mathbf{x}\in U\subset
\Sigma$ of $P$ where $X=\partial_t$, the flow $\phi^t$ acts by
translations in the $t$-direction. Thus, we may choose a
$\partial_t$-invariant frame $\{e_1,e_2,e_3\}$, $e_1=\partial_t=X$
on $V$ satisfying:
\[
[e_1,e_k]=[\partial_t,e_k]=0.
\]
Any local vector field $f$ on $U$ lifts, in a natural fashion, to the vector field $F$ on
$V\cong S^1\times U$, so that the equation $\pi_\ast(F)=f$ holds, and
we may define a metric $g_\Sigma$ on
$U\subset \Sigma$ by 
\[
g_\Sigma(f,f')=g_P(F,F').
\]
This turns $\pi$ into a Riemannian submersion on $V$, and defines $g_\Sigma$ pointwise on the whole $\Sigma$.
In the next step we obtain \eqref{eq:princ_sub2} as a corollary of Lemma \ref{th:bracket_submersion}.
Since the Christoffel symbols project under Riemannian submersions
(see \emph{(d)} in Lemma \ref{th:bracket_submersion}), for $u\in C^2(\Sigma)$ in
a local frame $\{E_2,E_3\}$ on $\Sigma$ we derive
\begin{eqnarray*}
(\Delta_\Sigma u)\circ \pi & = &(-D_{E_2}D_{E_2} u-D_{E_3}D_{E_3}
u+\Omega^3_{22}D_{E_3}
u+\Omega^2_{33}D_{E_2} u)\circ \pi\\
& = & -D_{e_2}D_{e_2}
(u\circ\pi)-D_{e_3}D_{e_3}(u\circ\pi)+\omega^3_{22}D_{e_3}
(u\circ\pi)+\omega^2_{33}D_{e_2}(u\circ\pi)\\
& = & \Delta_E (u\circ\pi),
\end{eqnarray*}
where $e_2=\mathfrak{H}(E_2)$ and $e_3=\mathfrak{H}(E_3)$.
\end{proof}

\no The following lemma (see \cite[p. 148]{Brin96}, Lemma 2.4.22) is
an important ingredient in the proof of the main theorem, thus we
provide a proof for more complete exposition.

\begin{lemma}[\cite{Brin96}]\label{th:seifert_fibered}
Every closed, compact orientable Seifert fibered 3-manifold $M$, with the base which is
a ``good'' orbifold $\Sigma$,
 is covered by a total space of a circle bundle
$P$. We have the following diagram:
\begin{align}\label{eq:seifert_diagram}
\CD
P @>p>> M\\
@VV{\Pi} V @VV{\pi}V\\
\tilde{\Sigma} @>r>> \Sigma
\endCD
\end{align}
where $p$ is the covering map, $r$ is the orbifold covering, and
the maps $\pi$, $\Pi$ are fibrations.
\end{lemma}
\begin{proof}
Any \emph{good} $2$-orbifold $\Sigma$ is a quotient of one of the model
spaces $S=S^2, \real^2$ or $\mathbb{H}^2$, by a discrete group of
isometries, i.e. $\Sigma=S/G$. Since $\Sigma$ is good then
any finitely generated discrete subgroup $G$ of isometries of $S$,
with compact quotient space, has a torsion free subgroup $G'$ of
finite index, \cite{Scott83}. Clearly, such a subgroup is isomorphic to the
fundamental group of a closed surface $\tilde{\Sigma}$. Define
$\tilde{\Sigma}=S/G'$, and $r:\tilde{\Sigma}\mapsto \Sigma$ to be a
quotient map, notice that $r$ is generally not a cover in the usual
sense, \cite{Scott83}. Let $h\in
\pi_1(M)$ represent a regular fiber of $M$. The subgroup $\langle
h\rangle$ of $\pi_1(M)$ generated by $h$ is infinite cyclic, and $\pi_1(M)/\langle
h\rangle=\pi_\ast(\pi_1(M))=G$. Denote by $K$ the inverse image in
$\pi_1(M)$ of a torsion free subgroup $G'$, under the induced group
homomorphism $\pi_\ast$. Let $P$ be the covering space of $M$
corresponding to $K$. If $\tilde{h}$ 
is represented in $\pi_1(P)=K\subset \pi_1(M)$ by a regular fiber then
$\langle\tilde{h}\rangle=\langle h\rangle$. Because $K/\langle
h\rangle=\pi_1(P)/\langle h\rangle=G'$ is torsion free, $P$ has no
singular fibers, and has to be an $S^1$-bundle over
$\tilde{\Sigma}$. Diagram \eqref{eq:seifert_diagram} follows
accordingly.
\end{proof}

\section{Proof of the Main Theorem.}\label{sec:6}
We now prove our main results. In a nutshell, these theorems describe
lower bounds for the volume of $M$ in terms of geometric parameters of $M$ and the magnitude $\|\alpha\|$ of a contact form
$\alpha$, which defines overtwisted or virtually overtwisted contact structure $\xi$ on $M$. Clearly, opposite inequalities provide sufficient conditions for the universal tightness of $\xi$. 
These results require existence of an $S^1$-action by a contact vector field which is 
also Killing in the adapted metric. Proposition \ref{th:seifert-adapted} demonstrates that, given an
$S^1$-invariant contact form $\alpha$, we may always adapt a
suitable Riemannian metric satisfying these requirements. Topologically, 
$M$ is a Seifert fibered manifold covered by a principal $S^1$-bundle $P$. Therefore, as
Theorem \ref{th:giroux_bundle} assures, the universal tightness of $\xi$ is
completely characterized by the topology of the dividing set on the base of $P$, and techniques developed in Section \ref{sec:4} can be applied. 

\begin{theorem}[Main Theorem (version 1)]\label{th:main}
\

\begin{itemize}
\item[\emph{(A)}] Let $(M,g_M)$ be a compact closed orientable Riemannian 3-manifold, equipped with a contact
structure $\xi$ defined by $\alpha$ and satisfying \eqref{eq:beltrami_form}. Assume that $\alpha$ admits a contact
vector field $X$ (i.e. $\Lie_X \alpha=0$) with circular orbits,
which is a unit Killing vector field for $g_M$.

\item[\emph{(B)}] Additionally, let the sectional curvature $\kappa_E$ of planes $E$, orthogonal
to the fibers, obey
\begin{align}\label{eq:negative-sectional}
\kappa_E\leq -\frac{3}{4}\,\mathcal{E}^2.
\end{align}
\end{itemize}

\no If $\xi$ is an overtwisted, or virtually overtwisted contact structure on $M$, then we have the following 
lower bound for the volume of $M$:

\begin{itemize}
\item[\emph{(C)}]
\begin{align}\label{eq:non-constant-tight}
\text{\rm Vol}(M)\geq \frac{2 \pi\, l_{\min}}{m_\alpha\, k},\qquad
m_\alpha=\max_M \bigl(0, \Delta_M\,\text{\emph{ln}}\|\alpha\|\bigr),
\end{align}
where $l_{\min}$ is a lower bound for lengths of orbits
of $X$, and $k$ depends on the Seifert fibration of $M$ induced by $X$.
\end{itemize}
\end{theorem}

\begin{proof}

Since $X$ has circular orbits, the result of Epstein \cite{Epstein72} implies that $M$ is a Seifert fibered manifold
and the lengths of orbits of $X$ are bounded. Consequently, $X$
induces an $S^1$-action by isometries on $M$, and we obtain an
orbifold bundle: $\pi:M\mapsto M/S^1 \simeq \Sigma$. 
In the first part of the proof, we show how (A) and (B) 
imply that the base $\Sigma$ of the Seifert fibration $M$ is a negatively curved orbifold. Such orbifolds are \emph{good} thus, $M$ is covered by an 
 $S^1$-bundle $P$, as concluded in Lemma \ref{th:seifert_fibered}. The constant $k$ is the degree of a cover. This allows us, in the second part of the proof, to lift the structure from 
 $M$ to the covering space $P$, using Diagram \eqref{eq:seifert_diagram}, and perform the analysis on $P$.

\bigskip

Let
$C=\{x_1,\ldots, x_k\}$ be the cone points of $\Sigma$, and
$S=\pi^{-1}(C)$ the set of singular fibers in $M$. Since $M\setminus
S\cong S^1\times (\Sigma\setminus C)$, by Lemma
\ref{th:princ_submersion} we may define a metric $g_\Sigma$ on
$\Sigma\setminus C$ so that $\pi:M\setminus S\mapsto \Sigma\setminus
C$ is a Riemannian submersion. The metric $g_\Sigma$ is smooth and
extends continuously to $\Sigma$. In the first step, we prove that the
scalar curvature of $(\Sigma\setminus C,g_\Sigma)$ is nonpositive, implying that
$\Sigma$ is a good orbifold (c.f. \cite{Scott83}).

Let us fix a local frame of vector fields $\{e_1=X,e_2,e_3\}$, and
the dual coframe $\{\eta=\eta_1,\eta_2,\eta_3\}$. Since $X$ is the
Killing vector field (i.e. $\Lie_X g_M=0$), for any pair of vector
fields $V$, $W$:
\begin{align*}
\langle D_V X, W\rangle=-\langle V, D_W X\rangle.
\end{align*}
Consequently, we obtain the following identities for the Christoffel
symbols in the frame $\{e_i\}$:
\begin{align}\label{eq:main_christ}
\omega^2_{1\,1}=\omega^3_{1\,1}=\omega^2_{2\,1}=\omega^3_{3\,1}=0,\qquad
\omega^k_{i\,j}=-\omega^j_{i\,k},\qquad 
 -\omega^2_{3\,1}=\omega^3_{2\,1}=\frac{\mathcal{E}}{2},
\end{align}
where $D_i e_j=\omega^k_{i\,j} e_k$.
Cartan's structure equations imply that the 1-form
$\eta=g_M(X,\,\cdot\,)$ satisfies
\begin{align}\label{eq:main_eta_eq}
\ast\,d\,\eta=\mathcal{E}\eta.
\end{align}
 Since $d\,\mathcal{E}\ast\eta=0$, we have $d\,\mathcal{E}(X)=X\,\mathcal{E}=0$,
 thus $\mathcal{E}$ is $S^1$-invariant. Using \eqref{eq:main_christ}, we compute
the sectional curvature $\kappa_E$ as follows (c.f.
\cite[p. 8]{Nicolaescu98})
\begin{eqnarray*}
\notag D_2 D_3 e_3 & = & D_2(\omega^k_{3\,3}\,e_k)=(D_2\,\omega^2_{3\,3})e_2+\omega^2_{3\,3}D_2 e_2\\
\notag & = &
(D_2\,\omega^2_{3\,3})e_2+\omega^2_{3\,3}\omega^3_{2\,2}e_3,\\
\notag D_3 D_2 e_3 & = & D_3(\omega^k_{2\,3}\,e_k)=(-\frac{1}{2}
D_3\,\mathcal{E})\,e_1
-\frac{1}{2}\mathcal{E} D_3 e_1+(D_3\,\omega^2_{2\,3})\,e_2+\omega^2_{2\,3}D_3 e_2\\
\notag & = & -\frac{1}{2} D_3\,\mathcal{E}\,
e_1+\frac{\mathcal{E}^2}{4} e_2
 +(D_3\,\omega^2_{2\,3}) e_2+\frac{\mathcal{E}}{2}\omega^3_{2\,2} e_1+\omega^3_{2\,2} \omega^2_{3\,3} e_3,
\end{eqnarray*}
\begin{eqnarray*}
\notag\left[e_2,e_3\right] & = & D_2 e_3-D_3 e_2= (\omega^k_{2\,3}-\omega^k_{3\,2})\,e_k=-\mathcal{E} e_1+\omega^2_{2\,3}e_2-\omega^3_{3\,2}e_3,\\
\notag D_{\left[e_2,\,e_3\right]} e_3 & = & -\mathcal{E} D_1
e_3+\omega^2_{2\,3}D_2
e_3-\omega^3_{3\,2}D_3 e_3\\
\notag & = & -\frac{\mathcal{E}^2}{2} \varphi
e_2+\frac{1}{2}\omega^2_{2\,3}\mathcal{E}
e_1+\left((\omega^3_{2\,2})^2 +(\omega^2_{3\,3})^2\right) e_2,
\end{eqnarray*}
\begin{eqnarray}\label{eq:sect_curv}
\kappa_E & = & \langle R(e_2,e_3)e_3,e_2\rangle\\
\notag & = & D_2\,\omega^2_{3\,3}+D_3\,\omega^2_{2\,3} -
\frac{\mathcal{E}^2}{4}+\mathcal{E}\varphi-\left((\omega^3_{2\,2})^2
+(\omega^3_{3\,2})^2\right)\\
 \notag & = & \sigma-\frac{\mathcal{E}^2}{4}+\mathcal{E}\varphi,
\end{eqnarray} 
where $\sigma=D_2\,\omega^2_{3\,3}+D_3\,\omega^2_{2\,3}-\left((\omega^3_{2\,2})^2
+(\omega^3_{3\,2})^2\right)$, and $\varphi=\omega^{2}_{1\,3}$. Notice that
\begin{align*}
D_1 e_2 = \omega^{3}_{1\,2}\, e_3=-\varphi\, e_3,\qquad D_1 e_3 =
\omega^{2}_{1\,3}\, e_2=\varphi\, e_2\,.
\end{align*}
Therefore, $\varphi$ measures a rotation of the frame in $E$, when
parallel transported along orbits of $X$. By Lemma
\ref{th:bracket_submersion}, the Christoffel symbols project under
$\pi: (M\setminus S,g_M)\mapsto (\Sigma\setminus C, g_\Sigma)$, and the scalar curvature $K$ of
$\Sigma$ obeys:
\[
K\circ\pi(x)=\sigma(x),\qquad \text{for}\quad x\in M,
\]
where $\sigma$ is defined in \eqref{eq:sect_curv}. Assuming that
$\{e_2, e_3\}$ are horizontal lifts of a frame from $\Sigma$,
Equation \eqref{eq:princ_sub1} yields
\[
0=[e_2, e_1]=D_1 e_2 - D_2 e_1,
\]
thus
\[
\varphi=-\frac{\mathcal{E}}{2}\,. 
\]
Because
$\kappa_E=\sigma-\frac{\mathcal{E}^2}{4}+\mathcal{E}\,\varphi=\sigma-\frac{3}{4}\,\mathcal{E}^2$,
by the assumption (B):
\begin{align}\label{eq:K_s}
K\circ\pi=\kappa_E+\frac{3}{4}\,\mathcal{E}^2 \leq 0.
\end{align}
By the Gauss-Bonnet theorem for orbifolds \cite{Scott83}: $\chi_\text{orb}(\Sigma)\leq 0$, and
$\Sigma$ must be covered by a closed surface $\tilde{\Sigma}$ of
nonzero genus (denote the covering projection by
$r:\tilde{\Sigma}\mapsto \Sigma$). Now, Lemma \ref{th:seifert_fibered}
tells us how to choose a principal bundle $\Pi: P\mapsto
\tilde{\Sigma}$, such that the total space $P$ is a covering space
for $M$. This is done as follows, recall that Diagram
\eqref{eq:seifert_diagram} commutes, and $p:P\mapsto M$ is a fiber
preserving covering map. Define a metric $g_P$ on $P$ by pulling
back the metric $g_M$ from $M$ via $p$, this makes $p:(P,g_P)\mapsto
(M,g_M)$ into a local isometry, and $\Pi:P\longrightarrow \tilde{\Sigma}$
into a Riemannian submersion. Let $\tilde{X}$ be the unique lift of
$X$, because $p$ respects the fibers, which are orbits of the flow
$\phi_X$ of $X$, we have
\begin{align*}
\phi_X(t,\cdot)\circ p=p\circ\phi_{\tilde{X}}(t,\cdot).
\end{align*}
Clearly, the lift $\tilde{X}$ of $X$ must also be a Killing
vector field on $P$ with circular orbits.

\bigskip

In the second part of the proof, we show the bound in (C), under the assumption that $\xi$ is 
overtwisted or virtually overwisted. 
Notice that it suffices to work with
the $S^1$-invariant contact structure $\tilde{\xi}$ on $P$, obtained by lifting $\xi$ to $P$ (i.e. $\tilde{\xi}=\ker{\tilde{\alpha}}$,
$\tilde{\alpha}=p_\ast\,\alpha$), since $\tilde{\xi}$ cannot be universally tight either. 
But $\tilde{\Sigma}\neq S^2$ thus, by the necessary and sufficient condition (b.1) in Theorem \ref{th:giroux_bundle}, $\tilde{\xi}$ satisfies 

\begin{itemize}
\item[($\ast$)]  The dividing set
$\Gamma_{\tilde{\Sigma}}$ on $\tilde{\Sigma}$, which is a projection of the characteristic surface $\Gamma_X$ under $\Pi$, contains a contractible closed curve.
\end{itemize}

\no Because $\Lie_{\tilde{X}}
\tilde{\alpha}=0$ and $\ast\,d\,\tilde{\alpha}=\mu\,\tilde{\alpha}$,
Theorem \ref{th:elliptic_contact} and Theorem
\ref{th:char_surface_geom} imply that $\Gamma_{\tilde{X}}=f^{-1}(0)$,
and $f=\alpha(X)\circ p$ is an $S^1$-invariant solution to 
Equation \eqref{eq:gener_lap}. By Lemma \ref{th:princ_submersion},
the following equation for $f$ holds on $\tilde{\Sigma}$:

\begin{equation}\label{eq:ellip_Sigma_eq}
\Delta_{\tilde{\Sigma}}\,f+\tilde{\mu}(\tilde{\mathcal{E}}-\tilde{\mu})\,f=0,
\end{equation}
where $\tilde{\mathcal{E}}=\mathcal{E}\circ p$, and $\tilde{\mathcal{\mu}}=\mathcal{\mu}\circ p$.
The function $f$ cannot be a trivial solution, for the following topological reason: $f\equiv 0$
 implies that $\tilde{\xi}$ is tangent everywhere to the $S^1$-fibers of $P$. But for $S^1$-invaraint $\tilde{\alpha}$, it would violate the nonintegrablity condition \eqref{eq:nonintegrable}. By Proposition \ref{th:necc_contact}, we must have $\tilde{\mu}(\tilde{\mu}-\tilde{\mathcal{E}})\geq 0$ on $M$, and unless $\tilde{\mu}=\tilde{\mathcal{E}}$ on $M$, $f$ cannot be a constant function.
(When $f$ is constant then $X$ is equal to the Reeb field of $X_\alpha$ and we arrive at 
Corrolary \ref{th:E=mu}.)

 If $f\neq \text{const}$, $f$ must change sign on
$\tilde{\Sigma}$, and the dividing set
$\Gamma_{\tilde{\Sigma}}=\Pi(\Gamma_{\tilde{X}})$ is
nonempty. (Notice that Theorem \ref{th:char_surface_geom}
implies that curves $\Gamma_\Sigma$ cannot have self-intersections.) By condition ($\ast$) one of the domains $\tilde{\Omega}\in
\tilde{\Sigma}\setminus \Gamma_{\tilde{\Sigma}}$ is a disc
$\tilde{\Omega}\cong D^2$. Notice that the function $\|\tilde{\alpha}\|$ is $\tilde{X}$-invariant (where $\tilde{\alpha}=p_\ast\,\alpha$), applying the technique of Proposition \ref{th:ellip_perp} we obtain
\begin{eqnarray}
\notag 2\pi & = & 2\pi\chi(\tilde{\Omega}) =  \int_{\tilde{\Omega}} K + \int_{\tilde{\Omega}}
\notag \Delta_{\tilde{\Sigma}}
 \text{\rm ln}\|\tilde{\alpha}\|\\
\label{eq:2pi_1} & \leq & \text{Vol}(\tilde{\Omega})\max_{\tilde{\Sigma}} \bigl(0, \Delta_{\tilde{\Sigma}} \text{\rm ln}\|\tilde{\alpha}\|\bigr),
\end{eqnarray}
because $K\leq 0$ (by \eqref{eq:K_s}).
In the next step, we bound the area of $\tilde{\Omega}$. Since
$r:\tilde{\Sigma}\setminus r^{-1}(C)\mapsto \Sigma\setminus C$ is a
$k$-sheeted cover and $\Sigma$ is a quotient of $\tilde{\Sigma}$ by
a discrete subgroup of isometries, we obtain
\[
\notag \text{Vol}(\tilde{\Omega})  \leq 
\text{Vol}(\tilde{\Sigma}) =  k\,\text{Vol}(\Sigma),
\]
and 
\begin{eqnarray}
\notag \text{Vol}(M) & = & \int_M \eta\wedge\pi^\ast\omega
  =  \int_\Sigma \int_{S^1}
\eta(X)\pi^\ast\omega\\
\label{eq:volume-of-M-est} & = & \int_{\Sigma} l(x)\omega\geq l_{\text{min}}\,
\text{Vol}(\Sigma),
\end{eqnarray}
where $l_{\text{min}}=\min_{x\in \Sigma}\ l(x)$ for the ``length of the fiber
function'' $l:\Sigma\mapsto \real$, and $\omega$ is the volume form on $\Sigma$. Bounds in  \eqref{eq:2pi_1} and \eqref{eq:volume-of-M-est} yield

\begin{equation}
\label{eq:non-const-cond} 2\pi \leq \frac{k}{l_\text{min}} \text{Vol}(M) \max_M \bigl(0, \Delta_M \text{\rm ln}\|\alpha\|\bigr).
\end{equation}
For $m_\alpha$ defined in (C) we obtain Inequality \eqref{eq:non-constant-tight}.
\end{proof}
\begin{corollary}
If $m_\alpha=0$, then $\xi$ is universally tight. 
\end{corollary}
\begin{corollary}\label{th:main-constant}
When $\mathcal{E}$ and $\mu$ are constant functions and $X$ has constant length $l$ orbits, Inequality \emph{(C)} simplifies to 
 \begin{itemize}
\item[\emph{(D)}] 
\[
\text{\rm Vol}(M)\geq \frac{4 \pi\, l}{\mu^2\,k}+\frac{2\pi\,e(M)\,l}{\mu},
\]
where $e(M)$ is the Euler number of the Seifert fibration of $M$ induced by $X$.
\end{itemize}

\end{corollary}
\begin{proof}
 The 1-form $\tilde{\eta}=\langle \tilde{X},\,\cdot\,\rangle$ can be regarded as a connection form on $P$ (since $\Lie_{\tilde{X}}\, \tilde{\eta}=0$, and $\ker\tilde{\eta}$ is orthogonal to the $S^1$-fibers).
 Therefore, we obtain the following relation between the function $\mathcal{E}$ and the Euler number of $P$ (c.f. \cite[p. 75]{McDuff-Salamon98}):
 \[
  e(P)=\frac{1}{2\pi}\int_{\tilde{\Sigma}} \tilde{\mathcal{E}}=\frac{k}{2\pi}\int_{\Sigma} \mathcal{E},
 \]
thus,
\begin{equation}\label{eq:euler-number}
 \mathcal{E}=\frac{2\pi\,e(P)}{\text{Vol}(\Sigma)\, k}=\frac{2\pi\,e(M)}{\text{Vol}(\Sigma)}\ .
\end{equation}
\no By \eqref{eq:K_s}, $\tilde{\Sigma}$ has 
nonpositive curvature. Because $\mathcal{E}=\tilde{\mathcal{E}}$ and $\mu=\tilde{\mu}$ are constant, Equation
\eqref{eq:ellip_Sigma_eq} is an eigenequation and Proposition
\ref{th:nodal_curves_cond} together with derivation \eqref{eq:volume-of-M-est} yield
\begin{eqnarray*}\label{eq:main_clnodal}
 4\,\pi & \leq &
 \mu(\mu-\mathcal{E})\text{Vol}(\tilde{\Omega})\\
  & = & \frac{k}{l} \mu(\mu-\mathcal{E})\text{Vol}(M),
\end{eqnarray*}
substituting \eqref{eq:euler-number} for $\mathcal{E}$ proves the claim.
\end{proof}

\begin{remark}\label{rem:no-circular-orbits}
 In Theorem \ref{th:main}, it may be possible to drop the assumption of circular orbits
 of $X$.
By the compactness of the group of isometries
of $(M, g_M)$, one easily
 shows that there exists a regular Killing vector field $X_\varepsilon$ arbitrarily close to $X$ (c.f. \cite{Nicolaescu98}).
 One may expect that Equation \eqref{eq:gener_lap} will hold for
 $f_\varepsilon =\alpha(X_\varepsilon)$ with possibly an error term.
 Consequently, one could imagine an approximation argument, with
$\varepsilon\to 0$, showing that the limit function $f$:
$f_\varepsilon\to f$, is a solution to \eqref{eq:gener_lap}.
\end{remark}

  The curvature assumption (B), in Theorem \ref{th:main}, is necessary to carry out the argument
 on the covering space $P$ of $M$, and also simplifies the bound in \eqref{eq:non-const-cond}. By assuming that the base of the 
 Seifert fibration $M$ is a good orbifold, we obtain

\begin{theorem}[Main Theorem (version 2)]\label{th:main-ver2}
\ 

\begin{itemize} 
\item[\emph{(E)}] Let $M$, $\alpha$, $X$, $g_M$ satisfy the assumption \emph{(A)}, and let 
$M/S^1 \simeq \Sigma$ be a good orbifold covered by a smooth surface of nonzero genus.
\end{itemize}
Then,
\begin{itemize} 
\item[\emph{(F)}] The form $\alpha$ defines a universally tight contact structure on $M$, provided that the volume of $M$ obeys
\[
\text{\rm Vol}(M)<\frac{2 \pi\, l_{\min}}{m_\alpha\, k},\qquad
m_\alpha=\max_M \bigl(0, (\Delta_M\,\text{\emph{ln}}\|\alpha\| + \kappa_E+\frac{3}{4}\,\mathcal{E}^2)\bigr),
\]
where $k$ is the degree of a cover (as in Theorem \ref{th:main}).
\end{itemize} 
\end{theorem}
\begin{proof}
We may apply analogous reasoning as in the second part of the proof of Theorem \ref{th:main}. The only difference is that $K=\kappa_E+\frac{3}{4}\,\mathcal{E}^2$ may be positive on $\tilde{\Sigma}$ so we must adjust the estimate in \eqref{eq:2pi_1} as follows
 \begin{eqnarray*}
\notag 2\pi & = & 2\pi\chi(\tilde{\Omega}) =  \int_{\tilde{\Omega}} K + \int_{\tilde{\Omega}}
\notag \Delta_{\tilde{\Sigma}}
 \text{\rm ln}\|\tilde{\alpha}\|\\
 & \leq & \text{Vol}(\tilde{\Omega})\max_{\tilde{\Sigma}} \bigl(0,\Delta_{\tilde{\Sigma}} \text{\rm ln}\|\tilde{\alpha}\|+\kappa_E+\frac{3}{4}\,\mathcal{E}^2\bigr).
\end{eqnarray*}
Now, the condition (F) may be derived analogously, as an opposite inequality.
\end{proof}

The results presented do not fully address the case when $\Sigma=M/S^1$ is homeomorphic to $S^2$, and is either 
a \emph{bad} orbifold or all covering spaces are $S^1$-bundles over $S^2$. 
In these cases, universal
tightness is determined by conditions (b.2) and (b.3) in Theorem \ref{th:giroux_bundle}. Corollary \ref{th:E=mu} provides only a partial answer here, and a different geometric condition is needed. Also, a drawback of lifting structures to the covering space is that one must know the degree $k$ of a cover, or its upper bound. Such upper bounds may be hard to obtain in full generality, but may be known for particular types of orbifolds (c.f. \cite{Scott83}). The author will address 
these issuess in the future work.

\section{Conclusions}\label{sec:7}
In concluding remarks, we want to point out several examples which demonstrate effectiveness of 
proposed theorems. 

The case when the contact field $X$ is the Reeb field $X_\alpha$, and thus $\xi$ is transverse to $S^1$-fibers, is captured by the following 
\begin{corollary}\label{th:E=mu}
If $\mu=\mathcal{E}$, under the assumption \emph{(E)}, $\xi$ is universally tight. Thus, 
 every regular Sasakian 3-manifold is universally tight (see \cite[p. 150]{Belgun03}). 
 The standard tight contact structure on $S^3$ is universally tight.
\end{corollary}

\no But, we also obtain more general result 

\begin{corollary}
If a contact form $\alpha$ is of constant length ($\|\alpha\|=\text{const}$), then under assumptions \emph{(A)} and \emph{(B)}, $\alpha$ defines a universally tight contact structure.
\end{corollary}

\no As a subsequent corollary, we obtain the well known result \cite{Kanda97} concerning
tight contact structures on a 3-torus $T^3\cong S^1\times S^1\times
S^1$.

\begin{corollary}
Every contact form $\alpha_n=\cos(n z) dx+\sin(n z) dy$, $n\in
\mathbb{Z}$, in standard coordinates on the 3-torus $T^3$, defines a universally tight contact structure on $T^3$.
\end{corollary}
\begin{proof}
The flat metric is adapted to $\alpha_n$ for all $n$ (c.f. \cite{Etnyre-Ghrist00a}) and the coordinate vector fields $\partial_x$, $\partial_y$ are unit Killing
vector fields with circular orbits. In addition, $\mathcal{E}=0$, $\|\alpha_n\|=1$ and $\kappa_E=0$, thus 
$m_{\alpha_n}=0$ in Theorem \ref{th:main-ver2}.
\end{proof}

To further demonstrate usefulness of proposed results,
 we look at the problem of \emph{geometric tightness} from a slightly different perspective. Rather then adapting a metric to a contact structure, we consider an arbitrary closed Riemannian 3-manifold $(M,g_M)$. The standard elliptic theory implies that the eigenproblem ($\mu=\text{const}$):
\begin{equation}\label{eq:curl-eigenproblem}
 \ast d\,\alpha=\mu\,\alpha,\qquad \mu\neq 0,
\end{equation}
admits infinitely many solutions, in particular we may diagonalize the operator $\ast d$ in the orthonormal basis of curl eigenfields: $\{\alpha_i\}$. If a solution $\alpha$ is a nonvanishing 1-form, i.e. $\|\alpha\|\neq 0$, the derivation in \eqref{eq:nonintegrability-derived} shows that the distribution 
$\xi=\ker\alpha$ defines a contact structure on $M$. Clearly, the Riemannian metric $g_M$ is adapted to every nonvanishing curl
eigenfield.

\begin{question}
 \emph{When are these contact structures tight/overtwisted}? 
\end{question}

If $(M,g_M)$ admits nonsingular Killing fields it also admits an $S^1$-action by the group of isometries (see Remark \ref{rem:no-circular-orbits}).
As shown in \cite{GK06}, $\ast d$ has the $S^1$-invariant portion of the spectrum, thus
\eqref{eq:curl-eigenproblem} admits $S^1$-invariant solutions. In the case of a unit Killing field,
Theorem \ref{th:main} and \ref{th:main-ver2} address the above question. 
To further demonstrate, let us consider a product of $S^1$ and a closed surface $\Sigma$: $M=S^1\times\Sigma$. Assume that $M$ has a product metric  $g_{S^1\times \Sigma}$ and $S^1$-fibers are of constant length $l$. In such setting we have a vertical vector field $X$ which is unit Killing in the metric $g_{S^1\times \Sigma}$. Since $X$ is orthogonal to $\Sigma$, the dual 1-form $\eta$ is closed and we obtain $\mathcal{E}=0$, by \eqref{eq:main_eta_eq}. If $\Sigma$ has nonpositive curvature, i.e. $K=\kappa_E\leq 0$, conditions (A) and (B) of Theorem \ref{th:main} are satisfied. Because $\mu$ and $\mathcal{E}$ are constant, every $S^1$-invariant solution to \eqref{eq:curl-eigenproblem} is universally tight provided 
\begin{equation}\label{eq:tight-condition-on_product}
 \text{Vol}(M)< \frac{4\pi l}{\mu^2},
\end{equation}
(compare to Corollary \ref{th:nodal_curves_cond_neg}). It can be easily shown
(see Lemma 4.2, \cite[p. 46]{GK06}) that eigenvalues $\mu^2$ of $S^1$-invariant curl eigenfields are equal to  eigenvalues $\lambda^2$ of the surface (in \eqref{eq:eigen-problem}). Therefore, Inequality \eqref{eq:tight-condition-on_product} will hold if $\Sigma$ has small eigenvalues with respect to the area of $\Sigma$.

Examples of Buser \cite{Buser77} show that hyperbolic surfaces with small eigenvalues are not uncommon. Specifically, Buser constructs a family of hyperbolic genus $g(\Sigma)$ surfaces (so called \emph{L$\ddot{o}$bell surfaces}) with $\text{Vol}(\Sigma)=4(g(\Sigma)-1)\pi$, where the eigenvalues satisfy
\[
 \lambda_i\leq \varepsilon,\qquad\text{for}\quad i=1,\,\ldots\,, 2 g(\Sigma)-3,
\]
(for arbitrary small $\varepsilon>0$). As a consequence, Inequality \eqref{eq:tight-condition-on_product} holds and the corresponding curl eigenfields are universally tight. In a nutshell, $\lambda_i$-eigenfunctions, for $i=1,\,\ldots\,, 2 g(\Sigma)-3$
cannot have a nodal domain which bounds a disc, and the corresponding curl eigenfields cannot have contractible dividing curves on $\Sigma$, which implies universal tightness. 
Consult \cite[p. 43]{Komendarczyk_thesis} for detailed considerations 
and a specific condition on \emph{L$\ddot{o}$bell surfaces} in Theorem 2.6.17. Techniques introduced in \cite{Komendarczyk_thesis, GK06} show that one may produce analogous examples in the setting of arbitrary $S^1$-bundles.

 These considerations confirm that knowledge of eigenvalues and the geometry of a manifold, in certain symmetric situations, is sufficient to determine  tightness of curl eigenfields without an explicit knowledge of solutions to the problem \eqref{eq:curl-eigenproblem}.  In general, this is not the case that for low eigenvalues 
nonsingular curl eigienfields define tight contact structures, \cite{GK06}. However, one may still hope that this claim \cite[p. 17]{Etnyre-Ghrist00} is true in a certain large class of Riemannian manifolds.

\section{Acknowledgments}
I am grateful to John Etnyre and Margaret Symington for discussions and support. I also wish to thank the  first referee for detailed corrections to the article.

\end{document}